\documentclass[12pt]{article}%
\usepackage{amsfonts}
\usepackage{sw20bams}
\usepackage{amsmath}
\usepackage{amssymb}
\usepackage{graphicx}%
\providecommand{\U}[1]{\protect\rule{.1in}{.1in}}
\begin{document}

\title{Number of Sign Changes: Segment of AR(1)}
\author{Steven Finch}
\date{September 5, 2019}
\maketitle

\begin{abstract}
Let $X_{t}$ denote a stationary first-order autoregressive process. \ Consider
$n$ contiguous observations (in time $t$) of the series (e.g., $X_{1}$,
\ldots, $X_{n}$). Let its mean be zero and its lag-one serial correlation be
$\rho$, which satisfies $\left\vert \rho\right\vert <1$. \ Rice (1945) proved
that $(n-1)\arccos(\rho)/\pi$ is the expected number of sign changes.
\ A\ corresponding formula for higher-order moments was proposed by Nyberg,
Lizana \& Ambj\"{o}rnsson (2018), based on an independent interval
approximation. \ We focus on the variance only, for small $n$, and see a
promising fit between theory and model.

\end{abstract}

\footnotetext{Copyright \copyright \ 2019 by Steven R. Finch. All rights
reserved.}Given%
\[%
\begin{array}
[c]{ccccc}%
X_{t}=\rho\,X_{t-1}+\sqrt{1-\rho^{2}}\cdot\varepsilon_{t}, &  & -\infty
<t<\infty, &  & |\rho|<1
\end{array}
\]
where $\varepsilon_{t}$ is $N(0,1)$ white noise, the segment $(X_{1}%
,\ldots,X_{n})$ is Gaussian with vector mean and covariance matrix%
\[
R=\left(
\begin{array}
[c]{ccccccc}%
1 & \rho & \rho^{2} & \ldots & \rho^{n-3} & \rho^{n-2} & \rho^{n-1}\\
\rho & 1 & \rho & \ldots &  & \rho^{n-3} & \rho^{n-2}\\
\rho^{2} & \rho & 1 &  & \ldots &  & \rho^{n-3}\\
\vdots & \vdots &  & \ddots & \ddots & \vdots & \vdots\\
\rho^{n-3} &  & \vdots & \ddots & 1 & \rho & \rho^{2}\\
\rho^{n-2} & \rho^{n-3} &  & \ldots & \rho & 1 & \rho\\
\rho^{n-1} & \rho^{n-2} & \rho^{n-3} & \ldots & \rho^{2} & \rho & 1
\end{array}
\right)  .
\]
In particular, all variances are one and the correlation between $X_{i}$ and
$X_{j}$ is $\rho^{\left\vert j-i\right\vert }$. \ Define $S_{n}=\#\{i:1\leq
i<n$ and $X_{i}X_{i+1}<0\}$, the number of sign changes, and%
\[
p_{e}(R)=\mathbb{P}\left\{  (-1)^{e_{1}}X_{1}<0\text{, \ }(-1)^{e_{2}}%
X_{2}<0\text{, \ldots, }(-1)^{e_{n-1}}X_{n-1}<0\text{ and }(-1)^{e_{n}}%
X_{n}<0\right\}
\]
for any vector $e$ of $n$ bits. \ It is well known that \cite{Gu1-heu,
Gu2-heu, Ow-heu, Tn-heu, GB-heu}
\begin{align*}
p_{11}  &  =\mathbb{P}\left\{  X_{1}>0\text{ and }X_{2}>0\right\} \\
&  =\frac{1}{4}+\frac{1}{2\pi}\arcsin\left(  \rho\right) \\
&  =\mathbb{P}\left\{  X_{1}<0\text{ and }X_{2}<0\right\}  =p_{00}%
\end{align*}
and
\begin{align*}
p_{111}  &  =\mathbb{P}\left\{  X_{1}>0\text{, }X_{2}>0\text{ and }%
X_{3}>0\right\} \\
&  =\frac{1}{8}+\frac{1}{4\pi}\left[  \arcsin\left(  \rho\right)
+\arcsin\left(  \rho^{2}\right)  +\arcsin\left(  \rho\right)  \right] \\
&  =\mathbb{P}\left\{  X_{1}<0\text{, }X_{2}<0\text{ and }X_{3}<0\right\}
=p_{000}.
\end{align*}
Because%
\[
\operatorname{Cov}\left(
\begin{array}
[c]{c}%
-X_{1}\\
X_{2}%
\end{array}
\right)  =\left(
\begin{array}
[c]{cc}%
1 & -\rho\\
-\rho & 1
\end{array}
\right)  ,
\]
we have%
\begin{align*}
p_{10}  &  =\mathbb{P}\left\{  X_{1}>0\text{ and }X_{2}<0\right\} \\
&  =\frac{1}{4}-\frac{1}{2\pi}\arcsin\left(  \rho\right) \\
&  =\mathbb{P}\left\{  X_{1}<0\text{ and }X_{2}>0\right\}  =p_{01}%
\end{align*}
and hence%
\[
\mathbb{E}\left(  S_{2}\right)  =2p_{10}=\frac{1}{2}-\frac{1}{\pi}%
\arcsin\left(  \rho\right)  =\frac{1}{\pi}\left[  \frac{\pi}{2}-\arcsin\left(
\rho\right)  \right]  =\frac{\arccos(\rho)}{\pi},
\]%
\begin{align*}
\mathbb{V}\left(  S_{2}\right)   &  =2p_{10}-\left(  2p_{10}\right)
^{2}=2p_{10}\left(  1-2p_{10}\right)  =\left[  \frac{1}{2}-\frac{1}{\pi
}\arcsin\left(  \rho\right)  \right]  \left[  \frac{1}{2}+\frac{1}{\pi}%
\arcsin\left(  \rho\right)  \right] \\
&  =\frac{1}{4}-\frac{1}{\pi^{2}}\arcsin\left(  \rho\right)  ^{2}.
\end{align*}
Because%
\[
\operatorname{Cov}\left(
\begin{array}
[c]{c}%
-X_{1}\\
X_{2}\\
X_{3}%
\end{array}
\right)  =\left(
\begin{array}
[c]{ccc}%
1 & -\rho & -\rho^{2}\\
-\rho & 1 & \rho\\
-\rho^{2} & \rho & 1
\end{array}
\right)
\]
and%
\[
\operatorname{Cov}\left(
\begin{array}
[c]{c}%
X_{1}\\
X_{2}\\
-X_{3}%
\end{array}
\right)  =\left(
\begin{array}
[c]{ccc}%
1 & \rho & -\rho^{2}\\
\rho & 1 & -\rho\\
-\rho^{2} & -\rho & 1
\end{array}
\right)
\]
we have%
\begin{align*}
p_{100}  &  =p_{011}=\frac{1}{8}+\frac{1}{4\pi}\left[  -\arcsin\left(
\rho\right)  -\arcsin\left(  \rho^{2}\right)  +\arcsin\left(  \rho\right)
\right] \\
&  =p_{001}=p_{110};
\end{align*}
because%
\[
\operatorname{Cov}\left(
\begin{array}
[c]{c}%
X_{1}\\
-X_{2}\\
X_{3}%
\end{array}
\right)  =\left(
\begin{array}
[c]{ccc}%
1 & -\rho & \rho^{2}\\
-\rho & 1 & -\rho\\
\rho^{2} & -\rho & 1
\end{array}
\right)
\]
we have%
\[
p_{010}=p_{101}=\frac{1}{8}+\frac{1}{4\pi}\left[  -\arcsin\left(  \rho\right)
+\arcsin\left(  \rho^{2}\right)  -\arcsin\left(  \rho\right)  \right]  ;
\]
thus%
\begin{align*}
\mathbb{E}\left(  S_{3}\right)   &  =4p_{100}+2\cdot2p_{010}\\
&  =\frac{1}{2}-\frac{1}{\pi}\arcsin\left(  \rho^{2}\right)  +\frac{1}%
{2}+\frac{1}{\pi}\left[  -2\arcsin\left(  \rho\right)  +\arcsin\left(
\rho^{2}\right)  \right] \\
&  =\frac{2}{\pi}\left[  \frac{\pi}{2}-\arcsin\left(  \rho\right)  \right]
=\frac{2\arccos(\rho)}{\pi},
\end{align*}%
\begin{align*}
\mathbb{V}\left(  S_{3}\right)   &  =4p_{100}+2\cdot4p_{010}-\left(
4p_{100}+2\cdot2p_{010}\right)  ^{2}\\
&  =\frac{1}{2}-\frac{1}{\pi}\arcsin\left(  \rho^{2}\right)  +1+\frac{2}{\pi
}\left[  -2\arcsin\left(  \rho\right)  +\arcsin\left(  \rho^{2}\right)
\right]  -\left[  1-\frac{2}{\pi}\arcsin\left(  \rho\right)  \right]  ^{2}\\
&  =\frac{3}{2}-\frac{4}{\pi}\arcsin\left(  \rho\right)  +\frac{1}{\pi}%
\arcsin\left(  \rho^{2}\right)  -1+\frac{4}{\pi}\arcsin\left(  \rho\right)
-\frac{4}{\pi^{2}}\arcsin\left(  \rho\right)  ^{2}\\
&  =\frac{1}{2}-\frac{4}{\pi^{2}}\arcsin\left(  \rho\right)  ^{2}+\frac{1}%
{\pi}\arcsin\left(  \rho^{2}\right)  .
\end{align*}
These formulas are consistent with a distributional result $S_{n}%
\sim\operatorname{Binomial}(n-1,1/2)$ valid when observations are independent;
in particular,
\[%
\begin{array}
[c]{ccc}%
\mathbb{E}\left(  S_{n}\right)  =\dfrac{n-1}{2}, &  & \mathbb{V}\left(
S_{n}\right)  =\dfrac{n-1}{4}%
\end{array}
\]
for $\rho=0$. \ The case $n=4$ for $\rho\neq0$ is more difficult and will be
covered in the next section. \ Closed-form variance expressions become
impossible for $n\geq5$ (see the appendix)\ and a certain approximative model
shall occupy us for the remainder of this paper.

\section{Dilogarithm Formula}

Cheng \cite{Ch1-heu, Ch2-heu, Ch3-heu, NK-heu} evaluated the following
integral:%
\[
I(h,x)=%
{\displaystyle\int\limits_{0}^{x}}
\arcsin\left(  \frac{\left(  1-h^{2}\right)  t}{h^{2}-t^{2}}\right)  \frac
{1}{\sqrt{1-t^{2}}}\,dt
\]
to be:%
\begin{align*}
&  -\frac{1}{2}\arcsin(x)^{2}+\frac{1}{2}\operatorname{Li}_{2}\left[
-\frac{\left(  h^{2}-\sqrt{h^{4}-x^{2}}\right)  ^{2}}{x^{2}}\right] \\
&  +\operatorname{Li}_{2}\left[  \frac{\left(  x-i\sqrt{1-x^{2}}\right)
\left(  h^{2}-\sqrt{h^{4}-x^{2}}\right)  }{x}\right]  +\operatorname{Li}%
_{2}\left[  \frac{\left(  x+i\sqrt{1-x^{2}}\right)  \left(  h^{2}-\sqrt
{h^{4}-x^{2}}\right)  }{x}\right] \\
&  -\frac{1}{2}\operatorname{Li}_{2}\left[  \frac{\left(  h^{2}-i\sqrt
{1-h^{2}}\right)  ^{2}\left(  h^{2}-\sqrt{h^{4}-x^{2}}\right)  ^{2}}{x^{2}%
}\right]  -\frac{1}{2}\operatorname{Li}_{2}\left[  \frac{\left(  h^{2}%
+i\sqrt{1-h^{2}}\right)  ^{2}\left(  h^{2}-\sqrt{h^{4}-x^{2}}\right)  ^{2}%
}{x^{2}}\right]
\end{align*}
where $0<x<h^{2}<1$ and $\operatorname{Li}_{2}[z]$ is the complex dilogarithm
function. Associated with covariance matrix%
\[
R^{+}=\left(
\begin{array}
[c]{cccc}%
1 & a & ab & a^{2}b\\
a & 1 & b & ab\\
ab & b & 1 & a\\
a^{2}b & ab & a & 1
\end{array}
\right)
\]
is orthant probability%
\[
p_{1111}(R^{+})=\frac{1}{16}+\frac{2\arcsin(a)+\arcsin(b)+2\arcsin
(ab)+\arcsin\left(  a^{2}b\right)  }{8\pi}+\frac{\arcsin(a)^{2}+I(a,a^{2}%
b)}{4\pi^{2}};
\]
call this $f(a,b)$. \ Associated with covariance matrix%
\[
R^{-}=\left(
\begin{array}
[c]{cccc}%
1 & -a & -ab & -a^{2}b\\
-a & 1 & b & ab\\
-ab & b & 1 & a\\
-a^{2}b & ab & a & 1
\end{array}
\right)
\]
is orthant probability%
\[
p_{1111}(R^{-})=\frac{1}{16}+\frac{\arcsin(b)-\arcsin\left(  a^{2}b\right)
}{8\pi}-\frac{\arcsin(a)^{2}+I(a,a^{2}b)}{4\pi^{2}};
\]
call this $g(a,b)$. \ We assume that $\left\vert a\right\vert <1$ and
$\left\vert b\right\vert <1$. \ Note that the matrix elements $R_{12}^{+}$ and
$R_{34}^{+}$ are identical, whereas $R_{12}^{-}$ and $R_{34}^{-}$ are of
opposite sign. \ Let us now return to our original $4\times4$ matrix $R$.
\ Clearly%
\begin{align*}
p_{1111}  &  =\mathbb{P}\left\{  X_{1}>0\text{, }X_{2}>0\text{, }X_{3}>0\text{
and }X_{4}>0\right\} \\
&  =f(\rho,\rho)\\
&  =\mathbb{P}\left\{  X_{1}<0\text{, }X_{2}<0\text{, }X_{3}<0\text{ and
}X_{4}<0\right\}  =p_{0000}.
\end{align*}
Because%
\[
\operatorname{Cov}\left(
\begin{array}
[c]{c}%
-X_{1}\\
X_{2}\\
X_{3}\\
X_{4}%
\end{array}
\right)  =\left(
\begin{array}
[c]{cccc}%
1 & -\rho & -\rho^{2} & -\rho^{3}\\
-\rho & 1 & \rho & \rho^{2}\\
-\rho^{2} & \rho & 1 & \rho\\
-\rho^{3} & \rho^{2} & \rho & 1
\end{array}
\right)
\]
and%
\[
\operatorname{Cov}\left(
\begin{array}
[c]{c}%
X_{1}\\
X_{2}\\
X_{3}\\
-X_{4}%
\end{array}
\right)  =\left(
\begin{array}
[c]{cccc}%
1 & \rho & \rho^{2} & -\rho^{3}\\
\rho & 1 & \rho & -\rho^{2}\\
\rho^{2} & \rho & 1 & -\rho\\
-\rho^{3} & -\rho^{2} & -\rho & 1
\end{array}
\right)
\]
we have%
\begin{align*}
p_{1000}  &  =p_{0111}=g(\rho,\rho)\\
&  =g(-\rho,\rho)=p_{0001}=p_{1110};
\end{align*}
because%
\[
\operatorname{Cov}\left(
\begin{array}
[c]{c}%
X_{1}\\
-X_{2}\\
X_{3}\\
X_{4}%
\end{array}
\right)  =\left(
\begin{array}
[c]{cccc}%
1 & -\rho & \rho^{2} & \rho^{3}\\
-\rho & 1 & -\rho & -\rho^{2}\\
\rho^{2} & -\rho & 1 & \rho\\
\rho^{3} & -\rho^{2} & \rho & 1
\end{array}
\right)
\]
and%
\[
\operatorname{Cov}\left(
\begin{array}
[c]{c}%
X_{1}\\
X_{2}\\
-X_{3}\\
X_{4}%
\end{array}
\right)  =\left(
\begin{array}
[c]{cccc}%
1 & \rho & -\rho^{2} & \rho^{3}\\
\rho & 1 & -\rho & \rho^{2}\\
-\rho^{2} & -\rho & 1 & -\rho\\
\rho^{3} & \rho^{2} & -\rho & 1
\end{array}
\right)
\]
we have
\begin{align*}
p_{0100}  &  =p_{1011}=g(\rho,-\rho)\\
&  =g(-\rho,-\rho)=p_{0010}=p_{1101}.
\end{align*}
Because%
\[
\operatorname{Cov}\left(
\begin{array}
[c]{c}%
-X_{1}\\
-X_{2}\\
X_{3}\\
X_{4}%
\end{array}
\right)  =\left(
\begin{array}
[c]{cccc}%
1 & \rho & -\rho^{2} & -\rho^{3}\\
\rho & 1 & -\rho & -\rho^{2}\\
-\rho^{2} & -\rho & 1 & \rho\\
-\rho^{3} & -\rho^{2} & \rho & 1
\end{array}
\right)
\]
we have
\[
p_{1100}=p_{0011}=f(\rho,-\rho);
\]
because%
\[
\operatorname{Cov}\left(
\begin{array}
[c]{c}%
X_{1}\\
-X_{2}\\
-X_{3}\\
X_{4}%
\end{array}
\right)  =\left(
\begin{array}
[c]{cccc}%
1 & -\rho & -\rho^{2} & \rho^{3}\\
-\rho & 1 & \rho & -\rho^{2}\\
-\rho^{2} & \rho & 1 & -\rho\\
\rho^{3} & -\rho^{2} & -\rho & 1
\end{array}
\right)
\]
we have
\[
p_{0110}=p_{1001}=f(-\rho,\rho);
\]
because%
\[
\operatorname{Cov}\left(
\begin{array}
[c]{c}%
-X_{1}\\
X_{2}\\
-X_{3}\\
X_{4}%
\end{array}
\right)  =\left(
\begin{array}
[c]{cccc}%
1 & -\rho & \rho^{2} & -\rho^{3}\\
-\rho & 1 & -\rho & \rho^{2}\\
\rho^{2} & -\rho & 1 & -\rho\\
-\rho^{3} & \rho^{2} & -\rho & 1
\end{array}
\right)
\]
we have
\[
p_{1010}=p_{0101}=f(-\rho,-\rho).
\]
Thus
\begin{align*}
\mathbb{E}\left(  S_{4}\right)   &  =4p_{1000}+2p_{1100}+4\cdot2p_{0100}%
+2\cdot2p_{0110}+2\cdot3p_{1010}\\
&  =4g(\rho,\rho)+2f(\rho,-\rho)+8g(\rho,-\rho)+4f(-\rho,\rho)+6f(-\rho
,-\rho)\\
&  =\frac{3\arccos(\rho)}{\pi},
\end{align*}%
\begin{align*}
\mathbb{V}\left(  S_{4}\right)   &  =4p_{1000}+2p_{1100}+4\cdot4p_{0100}%
+2\cdot4p_{0110}+2\cdot9p_{1010}-\left(  3\arccos(\rho)/\pi\right)  ^{2}\\
&  =4g(\rho,\rho)+2f(\rho,-\rho)+16g(\rho,-\rho)+8f(-\rho,\rho)+18f(-\rho
,-\rho)-9\left[  1/2-\arcsin(\rho)/\pi\right]  ^{2}.
\end{align*}
Unlike the mean, our expression for the variance does not simplify
appreciably. \ A plot of $\mathbb{V}\left(  S_{4}\right)  $ falls off
symmetrically from both sides of the maximum value $3/4$ at $\rho=0$. \ For
specificity's sake, we indicate numerical values at $\rho=1/2$:%
\[%
\begin{array}
[c]{ll}%
f(1/2,1/2)= & 0.1576625817544825416159596...,\\
g(1/2,1/2)= & 0.0707784073926423526601112...,\\
f(1/2,-1/2)= & 0.0658073315415406956707081...,\\
g(1/2,-1/2)= & 0.0390850126446677433865542...,\\
f(-1/2,1/2)= & 0.0341139367935660863971512...,\\
f(-1/2,-1/2)= & 0.0226893098357904842228499...,\\
\mathbb{V}\left(  S_{4}\right)  = & 0.7214075663610921033552384...<3/4
\end{array}
\]
and, of course,%
\[
2f(\rho,\rho)+4g(\rho,\rho)+2f(\rho,-\rho)+4g(\rho,-\rho)+2f(-\rho
,\rho)+2f(-\rho,-\rho)=1
\]
always.

\section{Independent Interval Approximation}

Our instinct (based on small samples) that the following should be true:%
\[%
\begin{array}
[c]{ccc}%
\mathbb{E}\left(  S_{n}\right)  =\dfrac{(n-1)\arccos(\rho)}{\pi} &  &
\text{for all }n\geq2
\end{array}
\]
is, in fact, a discrete-time analog of a classical theorem due to Rice
\cite{Ric-heu, Fi0-heu, Ked-heu, Bar-heu}.

The variance offers a more interesting situation. \ No pattern is evident from
our work and the case $n=5$ is beyond us. \ One tactic is to introduce a
modeling assumption that interval lengths between sign changes are
independently distributed. \ This idea apparently originated with Siegert
\cite{Sie-heu} and McFadden \cite{McF-heu} in the context of zero-crossings of
continuous-time processes, and suitably generalized in \cite{NAL-heu}. \ We
make no claim that the assumption is valid for most (or even some) processes.
\ It provides remarkably accurate estimates in many scenarios and our setting
is no exception.

Nyberg, Lizana \& Ambj\"{o}rnsson \cite{NLA-heu} obtained, within the
independent interval approximation (IIA) framework, a recursive formula%
\[%
\begin{array}
[c]{ccc}%
c_{n}=\dfrac{\arccos(\rho)}{6\pi}(n-1)n(n+1)-\dfrac{\pi}{\arccos(\rho)}%
{\displaystyle\sum\limits_{k=2}^{n-1}}
\left[  \dfrac{1}{2}-\dfrac{\arcsin\left(  \rho^{n-k+1}\right)  }{\pi}\right]
c_{k}, &  & c_{1}=0
\end{array}
\]
which is worthy of study. \ The quantity $c_{n}$ is the IIA-based estimate of
$\mathbb{E}\left(  S_{n}^{2}\right)  $. We calculate%
\[%
\begin{array}
[c]{ccc}%
c_{2}=\dfrac{\arccos(\rho)}{\pi}, &  & c_{3}=\dfrac{4\arccos(\rho)}{\pi
}-\dfrac{1}{2}+\dfrac{\arcsin\left(  \rho^{2}\right)  }{\pi}%
\end{array}
\]
and%
\[
c_{2}-\mathbb{E}\left(  S_{2}\right)  ^{2}=\dfrac{\arccos(\rho)}{\pi}\left[
1-\dfrac{\arccos(\rho)}{\pi}\right]  =\left[  \frac{1}{2}-\frac{\arcsin\left(
\rho\right)  }{\pi}\right]  \left[  \frac{1}{2}+\frac{\arcsin\left(
\rho\right)  }{\pi}\right]  =\mathbb{V}\left(  S_{2}\right)  ,
\]%
\begin{align*}
c_{3}-\mathbb{E}\left(  S_{3}\right)  ^{2}  &  =\dfrac{2\arccos(\rho)}{\pi
}\left[  2-\dfrac{2\arccos(\rho)}{\pi}\right]  -\dfrac{1}{2}+\dfrac
{\arcsin\left(  \rho^{2}\right)  }{\pi}\\
&  =-\dfrac{1}{2}+\left[  1-\frac{2\arcsin\left(  \rho\right)  }{\pi}\right]
\left[  1+\frac{2\arcsin\left(  \rho\right)  }{\pi}\right]  +\dfrac
{\arcsin\left(  \rho^{2}\right)  }{\pi}=\mathbb{V}\left(  S_{3}\right)  .
\end{align*}
That is, the model-based predictions of $\mathbb{V}\left(  S_{2}\right)  $ and
$\mathbb{V}\left(  S_{3}\right)  $ are exactly the same as theory! \ We also
calculate%
\[
c_{4}=\dfrac{10\arccos(\rho)}{\pi}-\dfrac{1}{2}+\dfrac{\arcsin\left(  \rho
^{3}\right)  }{\pi}-\dfrac{\pi}{\arccos(\rho)}\left[  \dfrac{1}{2}%
-\dfrac{\arcsin\left(  \rho^{2}\right)  }{\pi}\right]  \left[  \dfrac
{4\arccos(\rho)}{\pi}-\dfrac{1}{2}+\dfrac{\arcsin\left(  \rho^{2}\right)
}{\pi}\right]
\]
and here model $c_{4}-\mathbb{E}\left(  S_{4}\right)  ^{2}$ and theory
$\mathbb{V}\left(  S_{4}\right)  $ are not identical. \ The fit, however, is
promising (see Figure 1). \ The separation is largest ($\approx0.002$) for
positive $\rho$ when $\rho\approx0.763$; the separation is largest
($\approx0.036$) for negative $\rho$ when $\rho\approx-0.897$. \ 

The pronounced asymmetry in the model is inexplicable. \ We wonder if, in the
midst of elaborate IIA-based derivations, a positive correlation was
hypothesized (supported partly by the authors' decision \cite{NLA-heu} to
restrict their test simulations to $0<\rho<1$). \ Conceivably we are intended
to replace $\rho$ everywhere by $\left\vert \rho\right\vert $ in the formula
for $c_{n}$. \ This would force symmetry to occur and improve the fit. \ But
we are not certain of the intent.\footnote{Reasons underlying the hypothesis
$0<\rho<1$ may have to do more with historical context (in the physics
literature) than with any other factor.}

Higher-order moments were further discussed in \cite{NLA-heu}. \ The recursive
formula involving IIA-based estimates of $\mathbb{E}\left(  S_{n}^{3}\right)
$ is more complicated than that for $c_{n}$. \ It would be good someday to
implement this and to perform model-to-theory comparisons at the third-order
level, keeping the unresolved issue of negative correlation in mind.%
\begin{figure}[ptb]%
\centering
\includegraphics[
height=4.0127in,
width=6.0805in
]%
{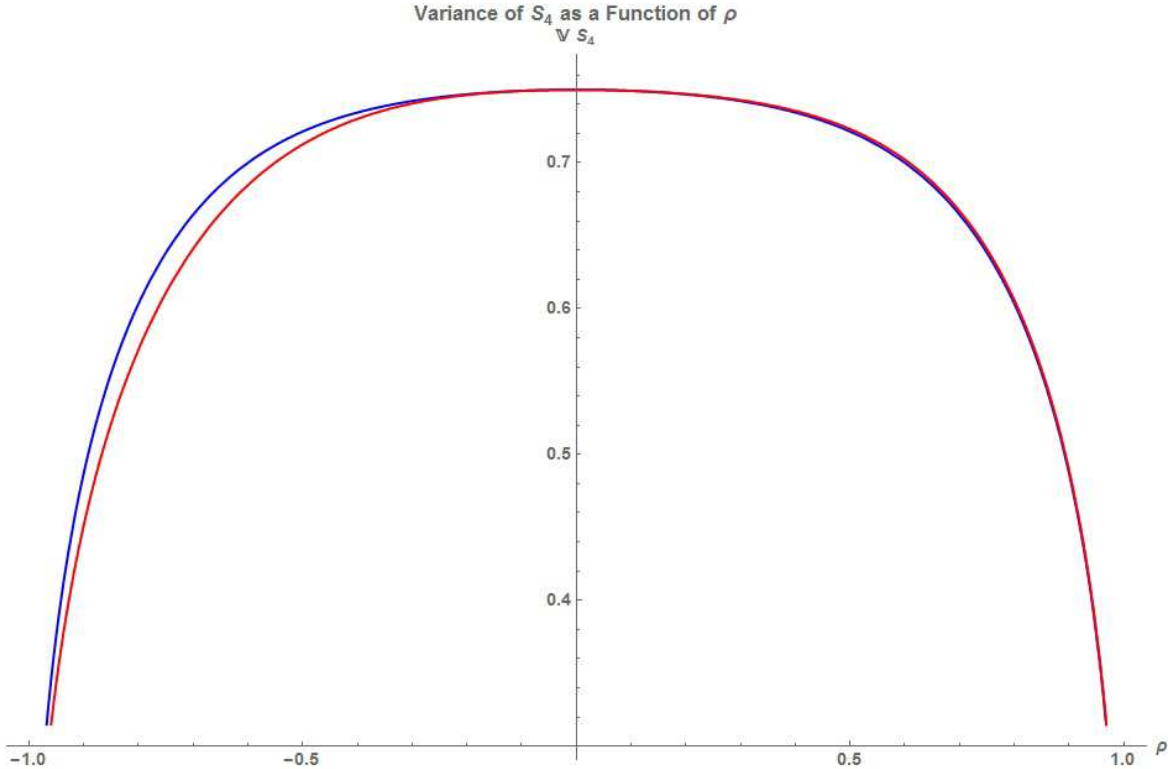}%
\caption{The red curve is an IIA-based model prediction of variance, while the
blue curve is our theoretical expression for variance. The blue curve is
symmetric with respect to the vertical axis; the red curve is not.}%
\end{figure}
\bigskip

\section{Appendix}

With regard to $n=5$, David \cite{Dvd-heu, Gui-heu} demonstrated how the
inclusion-exclusion principle can be applied to compute $p_{11111}=q_{12345}$.
\ Twenty-eight of the thirty terms in her expansion:%
\begin{align*}
q_{12345}  &  =\tfrac{1}{2}\left(  1-q_{1}-q_{2}-q_{3}-q_{4}-q_{5}%
+q_{12}+q_{13}+q_{14}+q_{15}+q_{23}+q_{24}+q_{25}+q_{34}+q_{35}+q_{45}\right.
\\
&  \;\;\;\;\;\;\;-q_{123}-q_{124}-q_{125}-q_{134}-q_{135}-q_{145}%
-q_{234}-q_{235}-q_{245}-q_{345}\\
&  \left.  \;\;\;\;\;\;\;+\,q_{1234}+q_{1235}+q_{1245}+q_{1345}+q_{2345}%
\right)
\end{align*}
can be easily evaluated. \ For example,%
\[
q_{1245}=\mathbb{P}\left\{  X_{1}>0\text{, }X_{2}>0\text{, }X_{4}>0\text{ and
}X_{5}>0\right\}
\]
possesses a closed-form expression because%
\[
\operatorname{Cov}\left(
\begin{array}
[c]{c}%
X_{1}\\
X_{2}\\
X_{4}\\
X_{5}%
\end{array}
\right)  =\left(
\begin{array}
[c]{cccc}%
1 & \rho & \rho^{3} & \rho^{4}\\
\rho & 1 & \rho^{2} & \rho^{3}\\
\rho^{3} & \rho^{2} & 1 & \rho\\
\rho^{4} & \rho^{3} & \rho & 1
\end{array}
\right)
\]
and this is of the form $R^{+}$ with $a=\rho$, $b=\rho^{2}$. \ The orthant
probability is%
\[
\frac{1}{16}+\frac{2\arcsin(\rho)+\arcsin\left(  \rho^{2}\right)
+2\arcsin\left(  \rho^{3}\right)  +\arcsin\left(  \rho^{4}\right)  }{8\pi
}+\frac{\arcsin(\rho)^{2}+I(\rho,\rho^{4})}{4\pi^{2}}%
\]
which is $0.1337768212694702494423619...$\ when $\rho=1/2$.

The two outlying terms:%
\[%
\begin{array}
[c]{ccc}%
\mathbb{P}\left\{  X_{1}>0\text{, }X_{2}>0\text{, }X_{3}>0\text{ and }%
X_{5}>0\right\}  , &  & \mathbb{P}\left\{  X_{1}>0\text{, }X_{3}>0\text{,
}X_{4}>0\text{ and }X_{5}>0\right\}
\end{array}
\]
are associated with matrices%
\[%
\begin{array}
[c]{ccc}%
\left(
\begin{array}
[c]{cccc}%
1 & \rho & \rho^{2} & \rho^{4}\\
\rho & 1 & \rho & \rho^{3}\\
\rho^{2} & \rho & 1 & \rho^{2}\\
\rho^{4} & \rho^{3} & \rho^{2} & 1
\end{array}
\right)  , &  & \left(
\begin{array}
[c]{cccc}%
1 & \rho^{2} & \rho^{3} & \rho^{4}\\
\rho^{2} & 1 & \rho & \rho^{2}\\
\rho^{3} & \rho & 1 & \rho\\
\rho^{4} & \rho^{2} & \rho & 1
\end{array}
\right)
\end{array}
\]
of a type so far unseen. \ The integral:%
\[
J(h,k,x)=%
{\displaystyle\int\limits_{0}^{x}}
\arcsin\left(  \frac{\sqrt{1-h^{2}}\sqrt{1-k^{2}}\,t}{\sqrt{h^{2}-t^{2}}%
\sqrt{k^{2}-t^{2}}}\right)  \frac{1}{\sqrt{1-t^{2}}}\,dt
\]
resists symbolic attack if $h\neq k$, but is nevertheless accessible to very
high-precision numerics. \ The two orthant probabilities are both equal to%
\[
\frac{1}{16}+\frac{2\arcsin(\rho)+2\arcsin\left(  \rho^{2}\right)
+\arcsin\left(  \rho^{3}\right)  +\arcsin\left(  \rho^{4}\right)  }{8\pi
}+\frac{\arcsin(\rho)\arcsin\left(  \rho^{2}\right)  +J(\rho,\rho^{2},\rho
^{4})}{4\pi^{2}}%
\]
which is $0.1354451520661386999235683...$\ when $\rho=1/2$. \ 

We close with two comments. \ First, our dilogarithm formula for $I(h,x)$
differs in appearance from Cheng's formula \cite{Ch1-heu} since he employed
$\operatorname{Li}_{2}[r,\theta]$ to represent the real part of
$\operatorname{Li}_{2}\left(  r\cdot e^{i\theta}\right)  $, whereas we use
\[%
\begin{array}
[c]{ccc}%
\operatorname{Li}_{2}\left(  z\right)  +\operatorname{Li}_{2}\left(
\overline{z}\right)  =2\operatorname{Re}\left[  \operatorname{Li}_{2}\left(
z\right)  \right]  &  & \text{for }z\in\mathbb{C}\setminus(1,\infty)
\end{array}
\]
to avoid this complication. \ Finally, given $0<k<1$, in an artificial
construct when%
\[%
\begin{array}
[c]{ccc}%
h=\dfrac{\sqrt{3}}{2}, &  & x=\sqrt{\dfrac{1-\sqrt{1-k^{2}}}{2}}%
\end{array}
\]
the integral $J(h,k,x)$ can be found \cite{Ch3-heu}:%
\begin{align*}
&  \frac{\pi^{2}}{8}-\frac{\pi}{6}\arcsin\left(  \ell\right)  +\frac{1}%
{6}\arcsin(\ell)^{2}-\frac{\pi}{2}\arcsin\left(  \sqrt{\frac{1-\ell}{2}%
}\right)  -\frac{1}{3}\operatorname{Li}_{2}\left[  -m^{2}\right] \\
&  -\frac{2}{3}\operatorname{Li}_{2}\left[  (\ell+ik)m\right]  -\frac{2}%
{3}\operatorname{Li}_{2}\left[  (\ell-ik)m\right]  +\frac{1}{3}%
\operatorname{Li}_{2}\left[  (\ell+ik)m^{2}\right]  +\frac{1}{3}%
\operatorname{Li}_{2}\left[  (\ell-ik)m^{2}\right]
\end{align*}
where%
\[%
\begin{array}
[c]{ccc}%
\ell=\sqrt{1-k^{2}}, &  & m=\dfrac{1+\ell-\sqrt{1+3\ell}\,\sqrt{1-\ell}}%
{2\ell}.
\end{array}
\]
This is a tantalizing hint that perhaps $J(\rho,\rho^{2},\rho^{4})$ is within
grasp if $\rho=\sqrt{3}/2$. \ Such a breakthrough will occur only if $x$ can
be unlocked from its current fixed location and allowed to wander free.

\section{Acknowledgements}

Helpful correspondence with Andreas Dieckmann \cite{Die-heu} and Tobias
Ambj\"{o}rnsson \cite{NLA-heu} is greatly appreciated. \ Virtually all
computations were performed within Mathematica. \ The pmvnorm function, part
of the mvtnorm package \cite{MH-heu} within R, was useful for verification
(impressive numerics: 12 digits of precision or better).

After the writing of this paper was completed, I\ learned of \cite{Sti-heu,
KS-heu}, which utilize similar techniques in answering somewhat different
questions. \ More aspects of AR(1) are covered in \cite{Fi1-heu, Fi2-heu}.

The cadence in much of Section 1 follows Emily Dickinson's lines
\textquotedblleft Because I could not stop for Death -- He kindly stopped for
me\textquotedblright. \ This paper is dedicated to the memory of my parents.

\end{document}